\documentclass[runningheads,a4paper]{llncs}

\usepackage[english]{babel}

\usepackage{amsfonts,amsmath}
\usepackage{amssymb}
\usepackage{eucal}
\usepackage{mathrsfs}
\usepackage[hypertexnames=false]{hyperref}
\usepackage{bm}
\usepackage{graphicx}
\usepackage{caption}
\usepackage{subcaption}
\usepackage{mathtools}
\usepackage{algorithm,algorithmic}

\usepackage[usenames,dvipsnames]{color}

\newcommand{\R}{\mathbb{R}}
\newcommand{\uu}{\hat{u}_k}
\DeclarePairedDelimiter{\ceil}{\lceil}{\rceil}

\providecommand{\norm}[1]{\left\lVert#1\right\rVert}

\begin{document}

\thispagestyle{empty}

\title{\mbox{Dynamic sampling schemes for optimal noise learning}\\
under multiple nonsmooth constraints}
\author{Luca Calatroni\inst{1}, Juan Carlos De Los Reyes\inst{2}, Carola-Bibiane Sch\"{o}nlieb\inst{3}}
\institute{Cambridge Centre for Analysis, University of Cambridge, UK \and Centro De Modelizaci\'{o}n Matem\'{a}tica, Escuela Polit\'{e}cnica Nacional de Quito, Ecuador \and Department of Applied Mathematics and Theoretical Physics, University of Cambridge, UK}

\titlerunning{Dynamic sampling schemes for optimal noise learning}
\authorrunning{L. Calatroni, J. C. De Los Reyes, C.-B. Sch\"{o}nlieb}
\date{}
\maketitle

\begin{abstract}
We consider the bilevel optimisation approach proposed in \cite{delosreyes} for learning the optimal parameters in a Total Variation (TV) denoising model featuring for multiple noise distributions. In applications, the use of databases (dictionaries) allows an accurate estimation of the parameters, but reflects in high computational costs due to the size of the databases and to the nonsmooth nature of the PDE constraints. To overcome this computational barrier we propose an optimisation algorithm that, by sampling \emph{dynamically} from the set of constraints and using a quasi-Newton method, solves the problem accurately and in an efficient way.
\end{abstract}
\vspace{-0.5cm}

\section{Introduction} \label{sec:int}

Most images in the real world suffer from noise. In photography noisy images occur when taking a photograph under bad lighting conditions, for instance. Medical imaging applications, such as Magnetic Resonance Imaging (MRI) and Positron Electron Tomography (PET), produce under-sampled and noisy image data. In general, the quality of images obtained from imaging devices in the real world, in the sciences and medicine, is limited by the hardware and the limited time available to measure the image data. Hence, one of the most important tasks in image processing is the reduction of noise in images, called image denoising.

A common challenge in image denoising is the setup of a suitable denoising model. The model depends on the noise distribution and the class of images the denoised solution should belong to.
%In many imaging applications such as Magnetic Resonance Imaging (MRI) or Positron Electron Tomography (PET) the accuracy, i.e. the noise level of the measurements, can be tuned by means of a database of real data. Such accuracy can depend, for instance, on the acquisition time of the machinery used. Thus, ground truth images or reasonably accurate approximations $u_k$ can be found by using the maximal accuracy possible, while their noisy versions $f_k$ can be measured within a usual clinical setup. Such a strategy is used already in the medical imaging community  where the good quality measurements or template shapes are used as priors in image reconstruction and image segmentation problems. The dataset can also be simulated: this is the case, for instance, of some MR velocity imaging applications where Boltzmann simulations need to be used to generate computer simulated data (for further details see \cite{mantle} and the references therein).
In \cite{delosreyes} a bilevel optimisation approach to learn the correct setup for a TV denoising model from a set of noisy and clean test images is proposed. There, optimal parameters $\lambda_i\in\R, i=1,\ldots,d$ are determined by solving the following optimisation problem:
\begin{equation}  \label{optprob1}
\min_{\lambda_i\geq 0,\  i=1,\ldots,d}\,\frac{1}{2N}\sum_{k=1}^{N}\norm{\uu-u_k}^2_{L^2(\Omega)}
%+\frac{\beta}{2}\sum_{i=1}^{d} |\lambda_i|^2
\end{equation}
subject to the set of nonsmooth constraints:
\begin{equation} \label{optprob2}
\uu=\text{argmin}_{u\in BV(\Omega)\cap\mathcal{A}}\left( |Du|(\Omega) + \sum_{i=1}^{d}\lambda_i\int_{\Omega} \phi_i(u,f_k)\,dx\right), \quad k=1,\ldots,N.
\end{equation}

In \eqref{optprob1}-\eqref{optprob2} $\Omega\subset\R^2$ is the image domain, $|Du|(\Omega)$ is the Total Variation (TV) $u$ in $\Omega$ and $BV(\Omega)$ is the space of functions of bounded variation (see \cite{AFP}). For each $k$, the pair $(u_k,f_k)$ is an element of a set of $N$ pairs of clean and noisy test images, respectively, whereas $\uu$ is the TV-denoised version of $f_k$. For $i=1,\ldots, d\ $ the terms $\phi_i$ represent the different data fidelities, each one modelling one particular type of noise weighted by a parameter $\lambda_i$. Examples of $\phi$ are $\phi(u,f_k)=(u-f_k)^2$ for noise with Gaussian distribution and $\phi(u,f_k)=|u-f_k|$ for the case of impulse noise. The set $\mathcal{A}$ is the set of admissible functions such that the data fidelity terms are well defined.

In this paper, we use a simulated database of clean and noisy images. This is not uncommon. Even in real world applications such as MRI, simulated databases are used to tune image retrieval systems, see for instance \cite{mantle}. Alternatively, we can imagine the retrieval of such a test set for a specific application using phantoms and their noisy acquisitions. Ideally, we would like to consider a very rich database (i.e. $N\gg 1$) in order to get a more robust estimation of the parameters, thus dealing with a very large set of constraints \eqref{optprob2} that would need to be solved in each iteration of an optimisation algorithm applied to  \eqref{optprob1}-\eqref{optprob2}. The computational solution of such an optimisation problem renders expensive and therefore challenging due to the large-scale nature of the problem \eqref{optprob1}-\eqref{optprob2} and due to the nonsmooth nature of each constraint.

%In order to deal with such large-scale problems various approaches have been presented in literature. They are based on the common idea of solving not \emph{all} the nonlinear PDE constraints, but just a sample of them, whose size varies according to the approach one intends to use. In \emph{stochastic approximation} (SA) methods (\cite{stochopt2,stochopt4}) typically a single data point is sampled per iteration, thus producing a generally inexpensive but noisy step. In \emph{sample} or \emph{batch average approximation} methods (see e.g. \cite{byrd2}) larger samples are used to compute an approximating (batch) gradient: the computation of such steps is normally more reliable, but generally more expensive. The  development of parallel computing, however, has improved upon the computational costs of batch methods: independent computation of functions and gradients can now be performed in parallel processors, so that the reliability of the approximation can be supported by more efficient methods.
%For optimisation problems under a large number constraints, some  efficient algorithms have been proposed in \cite{stochopt1,chung,stochopt3} under the name of Stochastic Average Approximation methods (SAA) and have shown good performance.
In order to deal with such large-scale problems various stochastic optimisation approaches have been presented in literature. They are based on the common idea of solving not \emph{all} the PDE constraints, but just a sample of them, whose size varies according to the approach one intends to use. In this paper we focus on a stochastic approximation method proposed by Byrd et al. \cite{byrd} called \emph{dynamic sample size} method. The main idea of this method is to consider an initial, small, training sample of the dictionary to start the algorithm with and \emph{dynamically} increasing its size, if needed, throughout the different steps of the optimisation process. The criterion to decide whether or not the sample size has to be increased  is a check on the sample variance estimates on the batch gradient. The desired trade-off between efficiency and accuracy is then obtained by starting with a small sample and gradually increasing its size till reaching the requested level of accuracy. Let us mention that the method of Byrd et al. is one among various stochastic optimisation methods, compare for instance \cite{stochopt2,stochopt4,stochopt1,chung,stochopt3}.
%the idea is to consider a Monte-Carlo approximation of the expected value of  the cost functional over all the random vectors $w$ with zero mean and covariance matrix $I$. Such approximation is found by considering $K$ different realisations of $w$ and $K$ is chosen such that $K\ll N$. Thus, this approach reduces each step of the optimisation algorithm to the solution of $K$ PDEs only and generally produces very accurate results as described in \cite{chung}. 

%In the following, we focus on another stochastic approximation technique, proposed by Byrd and al. in \cite{byrd} under the name of \emph{dynamic sample size} method. The main idea is to consider an initial, small, training sample of the dictionary to start the algorithm with and \emph{dynamically} increasing its size, if needed, throughout the different steps of the optimisation process. The criterion to decide whether or not the sample size has to be increased  is a check on the sample variance estimates on the batch gradient. The desired trade-off between efficiency and accuracy is then obtained by starting with a small sample and gradually increasing its size till reaching the requested level of accuracy. 

Our work extends the work of \cite{byrd} in two directions: firstly, in \cite{byrd} the linearity of the solution map is required which is not fulfilled for our problem \eqref{optprob1}-\eqref{optprob2}. We are going to show that the strategy of Byrd et al. can be modified for nonlinear solution maps as the one we are considering. Secondly, in \cite{byrd} the optimization algorithm is of gradient-descent type. Using a BFGS method to solve \eqref{optprob1}-\eqref{optprob2} we extend their approach incorporating also \emph{second order} information in form of an approximation of the Hessian by evaluations of the sample gradient in the iterations of the optimisation algorithm.

\medskip

\noindent{\textbf{Organisation of the paper.}} In the following Section \ref{sec:dynscheme} we present the Dynamic Sampling algorithm adapted to the nonlinear framework of problem \eqref{optprob1}-\eqref{optprob2}, specifying the variance condition on the batch gradient 
%and on the approximation of the Hessian
used in our optimisation algorithm. In Section \ref{sec:numres} we present the numerical results obtained for the estimation of the optimal parameters in the case of single and mixed noise estimation for the model \eqref{optprob1}-\eqref{optprob2} showing significant improvements in efficiency.
% We will show a significant improvement in the efficiency with respect to the algorithm solving the whole problem as well as a brief discussion regarding the sensitivity of the algorithm with respect to the accuracy parameter. 

\medskip

\noindent\textbf{Preliminaries.} We denote the vector of parameters we aim to estimate by $\bm{\lambda}=(\lambda_1,\ldots, \lambda_d)\in\R^d_{\geq 0}$. We also define by $\mathcal{S}$ the solution map that, for each constraint $k=1,\ldots,N$ of \eqref{optprob2}, associates to $\bm{\lambda}$ and to the noisy image $f_k$ the corresponding Total Variation denoised solution $\uu$, that is $\mathcal{S}(\bm{\lambda},f_k)=\uu$. Let us then define the reduced cost functional $J(\bm{\lambda})$ as
\begin{equation}   \label{functional}
J(\bm{\lambda}):=\frac{1}{2N}\sum_{k=1}^{N}\norm{\mathcal{S}(\bm{\lambda},f_k)-u_k}^2_{L^2(\Omega)}
%+\frac{\beta}{2}\sum_{i=1}^{d} |\lambda_i|^2.
\end{equation}
%As mentioned before, we assume that the size $N$ of the data set $\left\{(u_k,f_k)\right\}_{k=1,\ldots,N}$ is very large, which, combined to the nonlinear nature of the constraints in \eqref{optprob2}, makes the evaluation of $J(\bm{\lambda})$ extremely expensive.  
We also define:
\begin{equation} \label{lossfunction}
l(\bm{\lambda},f_k):=\norm{\mathcal{S}(\bm{\lambda},f_k)-u_k}^2_{L^2(\Omega)},\quad k=1,\dots,N
\end{equation}
as the \emph{loss functions} of the functional $J$ defined in \eqref{functional} for each $k=1,\ldots,N$. For every sample $S\subset\left\{1,\ldots,N\right\}$ of the database, we introduce the batch objective function:
\begin{equation} \label{batchfunction}
J_S(\bm{\lambda}):=\frac{1}{2|S|}\sum_{k\in S}l(\bm{\lambda},f_k)
%+\frac{\beta}{2}\sum_{i=1}^{d} |\lambda_i|^2.
\end{equation}

\section{Dynamic sampling schemes for solving \eqref{optprob1}-\eqref{optprob2}}  \label{sec:dynscheme}

To design the optimisation algorithm solving \eqref{optprob1}-\eqref{optprob2} we follow the approach used in \cite{delosreyes}. There, a quasi-Newton method (namely, the Broyden-Fletcher -Goldfarb-Shanno algorithm BFGS) is considered together with an Armijo backtracking linesearch rule. We combine such algorithm with a modified version of the Dynamic Sampling algorithm presented in \cite[Section 3]{byrd}.  In order to compare our algorithm with the Newton-Conjugate Gradient method presented in \cite[Section 5]{byrd}, we highlight that in our optimisation algorithm the Hessian matrix is never computed, but approximated efficiently by the BFGS matrix. 

Our algorithm starts by selecting from the whole dataset a sample $S$ whose size $|S|$ is small compared to the original size $N$. In the following iterations, if the approximation computed produces an improvement in the cost functional $J$, then the sample size is kept unchanged  and the optimisation process continues selecting in the next iteration a new sample of the same size. Otherwise, if the approximation computed is not a good one, a new, larger, sample size is selected and a new sample $S$ of this new size is used to compute the new step. By starting with small sample sizes it is hoped that in the early stages of the algorithm the solution can be computed efficiently in each iteration. The key point in this procedure is clearly the rule that checks throughout the progression of the algorithm, whether the approximation we are performing is good enough, i.e. the sample size is big enough, or has to be increased. Because of this systematic check on the quality of approximation in each step of the algorithm, such sampling strategy is called \emph{dynamic}.

As in \cite{byrd}, we consider a condition on the batch gradient $\nabla J_S$ which imposes in every stage of the optimisation that the direction $-\nabla J_S$ is a descent direction for $J$ at $\bm{\lambda}$ if the following condition holds:
\begin{equation} \label{descentcond}
\norm{\nabla J_S(\bm{\lambda})-\nabla J(\bm{\lambda})}_{L^2(\Omega)}
\leq \theta\norm{\nabla J_S(\bm{\lambda})}_{L^2(\Omega)},\quad \theta\in[0,1).
\end{equation}

The computation of $\nabla J$ may be very expensive for applications involving large databases and nonlinear constraints, so we rewrite \eqref{descentcond} as an estimate of the variance of the random vector $\nabla J_S(\bm{\lambda})$. In order to do that, recalling definitions \eqref{lossfunction} and \eqref{batchfunction} we observe that
\begin{equation}  \label{batchgradient}
\nabla J_S(\bm{\lambda})=\frac{1}{2|S|}\sum_{k\in S}\nabla l(\bm{\lambda},f_k)
%+\beta\sum_{i=1}^{d} \lambda_i.
\end{equation}
We can compute \eqref{batchgradient} in correspondence to an optimal solution $\bm{\hat{\lambda}}$ by using \cite[Remark 3.4]{delosreyes} where a characterisation of $\nabla J$ is given in terms of the adjoint states $p_k$ (see Section 3 for details). By linearity and extending to the multiple-constrained case, we get:
\begin{equation}   \label{batchgradientadj}
\nabla J_S(\bm{\hat{\lambda}})= \sum_{k\in S}\sum_{i=1}^{d}\int_\Omega \phi_i'(\uu,f_k) \,p_k\, dx
%+\beta\bm{\hat{\lambda}}.
\end{equation}
%\red{We refer the reader to Section \ref{sec:numres} where the computation of the adjoint state will be specified for each choice of $\phi_i$.}
Thanks to this characterisation, we now extend the dynamic sampling algorithm in \cite{byrd} to the case where the solution map $\mathcal{S}$ is nonlinear: by taking \eqref{batchgradientadj} into account and following \cite[Section 3]{byrd} we can rewrite \eqref{descentcond} as a condition on the variance of the batch gradient that reads as
\begin{equation} \label{descentdirectionvar}
\frac{\norm{Var_{k\in S}(\nabla l(\bm{\lambda},f_k))}_{L^1(\Omega)}}{|S|}\frac{N-|S|}{N-1}\leq\theta^2\norm{\nabla J_S(\bm{\lambda})}^2_{L^2(\Omega)}.
\end{equation}
For a detailed derivation of \eqref{descentdirectionvar}, see \cite{byrd}. Condition \eqref{descentdirectionvar}  is the responsible for possible changes in the sample size in the optimisation algorithm and has to be checked in every iteration. If inequality \eqref{descentdirectionvar} is not satisfied, a larger sample $\hat{S}$  whose size satisfies the descent condition \eqref{descentdirectionvar} needs to be considered. By assuming that the change in the sample size is gradual enough such that, for any given $\bm{\lambda}$:
\begin{align*}
&\norm{Var_{k\in \hat{S}}(\nabla l(\bm{\lambda},f_k))}_{L^1(\Omega)}\approx
\norm{Var_{k\in S}(\nabla l(\bm{\lambda},f_k))}_{L^1(\Omega)} ,\\
&\norm{\nabla J_{\hat{S}}(\bm{\lambda})}_{L^2(\Omega)}\approx \norm{\nabla J_S(\bm{\lambda})}_{L^2(\Omega)},
\end{align*}
we see that condition \eqref{descentdirectionvar} is satisfied whenever we choose $|\hat{S}|$ such that
\begin{equation} \label{samplesize}
|\hat{S}|\geq \ceil[\Bigg]{\frac{N-\norm{Var_{k\in S}(\nabla l(\bm{\lambda},f_k))}_{L^1(\Omega)}}{\norm{Var_{k\in S}(\nabla l(\bm{\lambda},f_k))}_{L^1(\Omega)}+\theta^2(N-1)\norm{\nabla J_S(\bm{\lambda})}^2_{L^2(\Omega)}}}.
\end{equation}
Conditions \eqref{descentdirectionvar} and \eqref{samplesize} are the key points in the optimisation algorithm we are going to present: by checking the former, one can control whether the sampling approximation is accurate enough and if this is not the case at any stage of the algorithm, by imposing the latter a new larger sample size is determined. 

\smallskip

We remark that these two conditions force the direction $-\nabla J_S$ to be a descent direction. Steepest descent methods are known to be slowly convergent. Algorithms incorporating information coming from the Hessian are generally more efficient. However, normally the computation of the Hessian is very expensive, so Hessian-approximating methods are commonly used. In \cite{byrd} a Newton-CG method is employed. There, an approximation of the Hessian matrix $\nabla^2 J_S$ is computed only on a subsample $H$ of $S$ such that $|H|\ll |S|$. As the sample $S$ is dynamically changing, the subsample $H$ will change as well (with a fixed, constant ratio) and the computation of the new conjugate gradient direction can be performed efficiently. In this work, in order to compute an approximation of the Hessian we consider the well-known BFGS method which has been extensively used in the last years because of its efficiency and low computational costs.
% \red{We stress here that it would not make sense considering our dynamic sampling approach in the BFGS optimisation in the unconstrained case. This because the sample considered and, consequently, the BFGS matrix, is changing in every iteration of BFGS. This is not a problem in our case though, as the dynamic sampling strategy is used \emph{just} for solving the PDE constraints: the BFGS is built after solving the PDEs thanks to the sampling strategy, so, when evaluating the quality of the approximation, comparing two different iterations still makes sense REPHRASE}.

Before giving a full description of the resulting algorithm solving \eqref{optprob1}-\eqref{optprob2}, we briefly comment on the linesearch rule that is employed in the update of the BFGS matrix.  We choose an Armijo backtracking line search rule with curvature verification: the BFGS matrix is updated only if the curvature condition is satisfied. The Amijo criterion is:
\begin{equation}   \label{armijo}
J_S(\bm{\lambda_k}+\alpha_k d_k)-J_S(\bm{\lambda_k})\leq \alpha_k\eta\nabla J_S(\bm{\lambda_k})^\top d_k
\end{equation}
where the value $\eta$ will be specified in Section \ref{sec:numres}, $d_k$ is the descent direction of the quasi-Newton step, $\alpha_k$ is the length of the quasi-Newton step and $\nabla J_S(\bm{\lambda_k})$ is defined in \eqref{batchgradient} . The positivity of the parameters is always preserved along the iterations.

\smallskip

We present now the BFGS optimisation with Dynamic sampling for solving \eqref{optprob1}-\eqref{optprob2}: compared to \cite[Algorithm 5.2]{byrd} we stress once more that the gain in efficiency is obtained thanks to the use of BFGS instead of the Newton-CG sampling method. 

\begin{algorithm}
\caption{Dynamic Sampling BFGS for solving \eqref{optprob1}-\eqref{optprob2}}
\begin{algorithmic}[1]
\STATE\small Initialize: $\bm{\lambda_0}$, sample $\mathcal{S}_0$ with $|S_0|\ll N$ and model parameter $\theta$, $k=0$.
\STATE \textbf{while} {BFGS not converging, $k\geq0$}
\STATE{\quad  sample $|S_k|$ PDE constraints to solve};
\STATE{\quad update of the BFGS matrix};
\STATE{\quad compute direction $d_k$ by BFGS and steplength $\alpha_k$ by Armijo cond. \eqref{armijo};}
\STATE{\quad define new iterate: $\bm{\lambda_{k+1}}=\bm{\lambda_k}+\alpha_k d_k$;}
\STATE{\quad \textbf{if} condition \eqref{descentdirectionvar} \textbf{then}}
\STATE{\quad\quad maintain the sample size: $|S_{k+1}|=|S_k|$;}
\STATE{\quad\textbf{else} augment $S_k$ such that condition \eqref{samplesize} is verified.}
\STATE{\textbf{end}}
\end{algorithmic}
\end{algorithm}

%\vspace{-1cm}

\section{Numerical results}  \label{sec:numres}
\setcounter{equation}{0}

In this section we present the numerical results of the Dynamic Sampling Algorithm 1 applied to compute the numerical solution of \eqref{optprob1}-\eqref{optprob2}. In our numerical computations we fix the parameter values as follows:
\begin{itemize}
\item We consider images of size $150\times 150$. We approximate the differential operators by discretising with finite difference schemes with mesh step size $h=1/$(number of pixels in the $x$-direction). We use forward difference for the discretisation of the divergence operator and backward differences for the gradient. The Laplace operator is discretised by using the usual five point formula.
\item The TV constraints in \eqref{optprob2} are solved by means of SemiSmooth Newton (SSN) algorithms whose form depends on the $\phi$'s in \eqref{optprob2} (cf. \cite[Section 4]{delosreyes}) solving regularised problems which stop if either the difference between two consecutive iterates is small enough or if the maximum number of iterations $\textsf{maxiter}=35$ is reached.
\item In the Armijo condition \eqref{armijo} the value $\eta$ is chosen to be $\gamma=10^{-4}$.
\end{itemize}

\paragraph{\textbf{Single noise estimation.}} As a toy example, we start by considering the case when the noise in the images is normally distributed. In \eqref{optprob1}-\eqref{optprob2}, this reflects in the estimation of just one parameter $\lambda$ that weights the fidelity term $\phi(u,f_k)=(u-f_k)^2$ in each constraint. Considering the training database $\left\{(u_k,f_k)\right\}_{k=1,\ldots,N}$ of clean and noisy images, the problem reduces to:
\begin{equation}  \label{optprob1gauss}
\min_{\lambda \geq 0}\,\frac{1}{2N}\sum_{k=1}^{N}\norm{\uu-u_k}^2_{L^2(\Omega)}
%+\frac{\beta}{2}\lambda^2
\end{equation}
where, for each $k$, $\uu$ is the solution of the regularised PDE
%a Huber-regularised version of the following elliptic-regularised variational inequalities 
%\begin{align} \label{optprob2gaussreg}
%-\varepsilon\Delta u_k^\gamma - \text{div}\Big{(}\frac{\gamma\nabla u^\gamma_k}{\text{max}(\gamma|\nabla u^\gamma_k|,1)}\Big{)} +\lambda (u_k^\gamma-f_k) =0 
%\end{align}
%where the second term relates to the derivative of a Huber-type regularisation of the subdifferential of $|Du^\gamma_k|$ with parameter $\gamma\gg 1$, for every $k$. The $\varepsilon$-term is a a small artificial diffusion term that sets up the problem in the Hilbert space $H^1_0(\Omega)$ (see \cite[Section 3]{delosreyes} for the theoretical details).}
%\blue{I would rather write (do you agree?)
\begin{align} \label{optprob2gaussreg}
-\varepsilon\Delta \uu - \text{div}\Big{(} h_\gamma( \nabla \uu) \Big{)} +\lambda (\uu-f_k) =0,\quad k=1,\dots,N.
\end{align}
In \eqref{optprob2gaussreg} $h_\gamma$ arises from a Huber-type regularisation of the subdifferential of $|D \uu|$ with parameter $\gamma\gg 1$ and the $\varepsilon$ term is an artificial diffusion term that sets up the problem in the Hilbert space $H^1_0(\Omega)$ (see \cite[Section 3]{delosreyes} for details). 

As shown in \cite[Theorem 3.5]{delosreyes} the adjoint states $p_k$ can be computed for each constraint as the solution of the following equation
\begin{multline} \label{eq:optimality condition adjoint equation}
 \varepsilon (D p_k, Dv)_{L^2}+ (h_\gamma'(D \uu)^* D p_k, D v)_{L^2}\\ + \int_\Omega \lambda ~p_k ~v ~dx =-(\uu-f_k,v)_{L^2},\quad \forall v \in H_0^1(\Omega).
\end{multline}
%\sum_{i=1}^d ~\phi_i''(\uu,f_k)
Recalling also equations \eqref{batchgradient}-\eqref{batchgradientadj} needed for the computation of the gradient, we can now apply Algorithm 1 to solve \eqref{optprob1gauss}-\eqref{optprob2gaussreg}.

For the following numerical tests, the parameters of this model are chosen as follows: $\varepsilon=10^{-12}, \gamma=100$. The noise in the images has distribution $\mathcal{N}(0,0.05)$. The parameter $\theta$ of the Algorithm 1, is chosen to be $\theta=0.5$. We will comment on the sensitivity of the method to $\theta$ later on.

\smallskip

Table \ref{tablegauss} shows the numerical value of the optimal parameter $\hat{\lambda}$ when varying the size of the dictionary. We measure the efficiency of the algorithms used in terms of the number of nonlinear PDEs solved during the BFGS optimisation and we compare the efficiency of solving \eqref{optprob1gauss}-\eqref{optprob2gaussreg} without and with the Dynamic Sampling strategy. We observe a clear improvement in efficiency when using Dynamic Sampling: the number of PDEs solved in the optimisation process is very much reduced. We note that this corresponds to an increasing number of BFGS iterations which does not appear to be an issue as BFGS iterations are themselves very fast. For the sake of computational efficiency, what really matters is the number of PDEs that need to be solved in \emph{each} iteration of BFGS. Moreover, thanks to modern parallel computing methods and to the decoupled nature of the constraints in each BFGS iteration, solving such a reduced amount of PDEs makes the computational efforts very reasonable. In fact, we note that the size of the sample is generally maintained very small in comparison to $N$ or just slightly increased. Computing also the relative error between the optimal parameter computed by solving all the PDEs and the one computed with Dynamic Sampling method, we note a good level of accuracy: the difference between the two values remains below $5\%$.
%\vspace{-0.5cm}
\begin{center}
{\scriptsize    
\begin{tabular}{| c | c | c | c | c | c | c | c | c | c | } 
    \hline
$N$ & $\hat{\lambda}$ & $\hat{\lambda}_S$ & $|S_0|$ & $|S_{end}|$ & eff. & eff. Dyn.S. & BFGS its. & BFGS its. Dyn.S.& diff. \\ \hline
$10$ & $3334.5$ & $3427.7$ &$2$&$3$& $140$ &\bm{$84$}& $7$ &$21$& $2.7\%$ \\ \hline
$20$ & $3437.0$ & $3475.1$ &$4$&$4$& $240$ &\bm{$120$}& $7$ &$15$& $1.1\%$ \\ \hline
$30$ & $3436.5$ & $3478.2$ &$6$&$6$& $420$ &\bm{$180$}& $7$ & $15$ & $1.2\%$ \\ \hline
$40$ & $3431.5$ & $3358.3$ &$8$&$9$& $560$ &\bm{$272$}& $7$ &$16$& $2.1\%$ \\ \hline
$50$ & $3425.8$ & $3306.4$ &$10$ & $10$ & $700$ & \bm{$220$} & $7$ & $11$ & $3.5\%$ \\ \hline
$60$ & $3426.0$ & $3543.4$ & $12$ & $12$ & $840$ & \bm{$264$} & $7$ &$11$ & $3.3\%$ \\ \hline
$70$ & $3419.7$ & $3457.7$ & $14$ & $14$ &  $980$ & \bm{$336$ }& $7$ & $12$ & $1.1\%$ \\ \hline
$80$ & $3418.1$ & $3379.3$  & $16$ & $16$ & $1120$ & \bm{$480$} &$7$ & $15$ & $<1\%$ \\ \hline
$90$ & $3416.6$ & $3353.5$  & $18$  & $18$  & $1260$ & \bm{$648$}  & $7$ & $18$  & $2.3\%$ \\ \hline
$100$ & $3413.6$ & $3479.0$ & $20$ & $20$ & $1400$ & \bm{$520$} & $7$ & $13$ & $1.9\%$ \\ \hline
 \end{tabular}}
\captionof{table}{$N$ is the size of the database, $\hat{\lambda}$ is the optimal parameter for \eqref{optprob1gauss}-\eqref{optprob2gaussreg} obtained by solving all the $N$ constraints, whereas $\hat{\lambda}_S$ is the one computed by solving the problem with Algorithm 1. The initial size $S_0$ is chosen to be $|S_0|=20\% N$. $|S_{end}|$ of the sample at the end of the optimisation algorithm. The efficiency of the  algorithms is measured in terms of the PDEs solved. We compare the accuracy of the  approximation in terms of the  difference $\norm{\hat{\lambda}_S-\hat{\lambda}}_1/\norm{\lambda_S|}_1$.}
\label{tablegauss}
\end{center}

\begin{figure}[h!]
\begin{center}
\includegraphics[width=2cm,height=2cm]{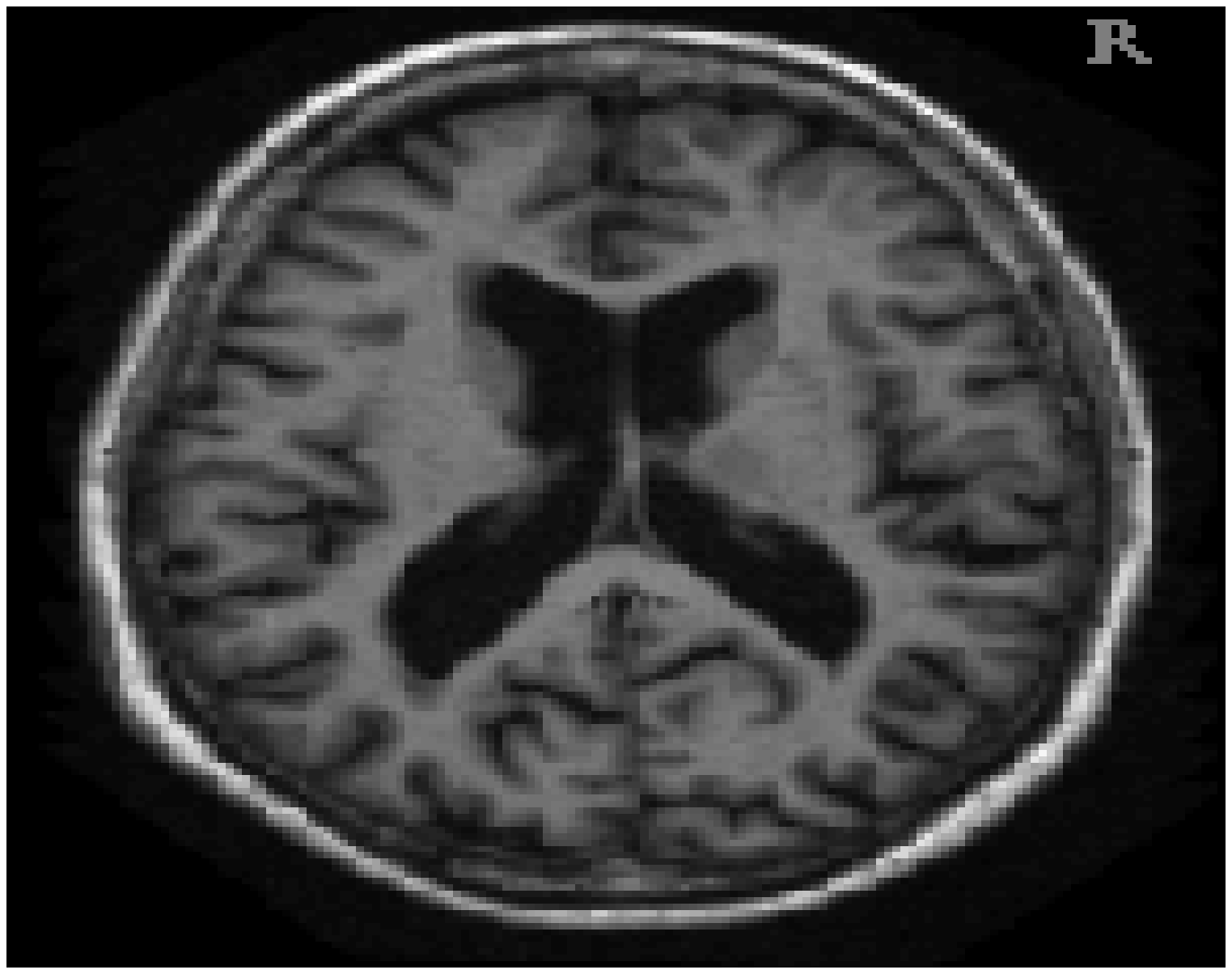}
\hspace{-0.3cm}
\includegraphics[width=2cm,height=2cm]{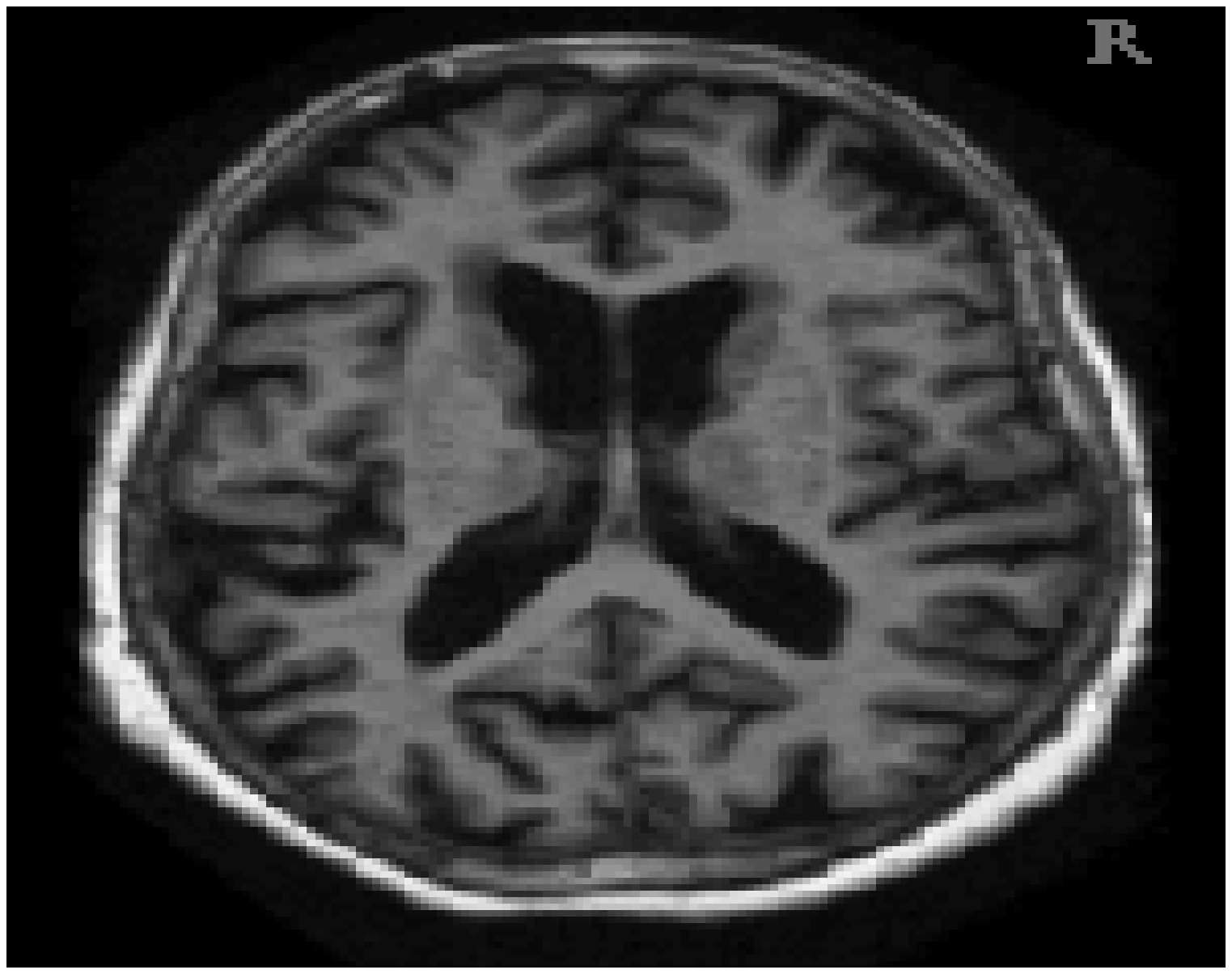}
\hspace{-0.3cm}
 \includegraphics[width=2cm,height=2cm]{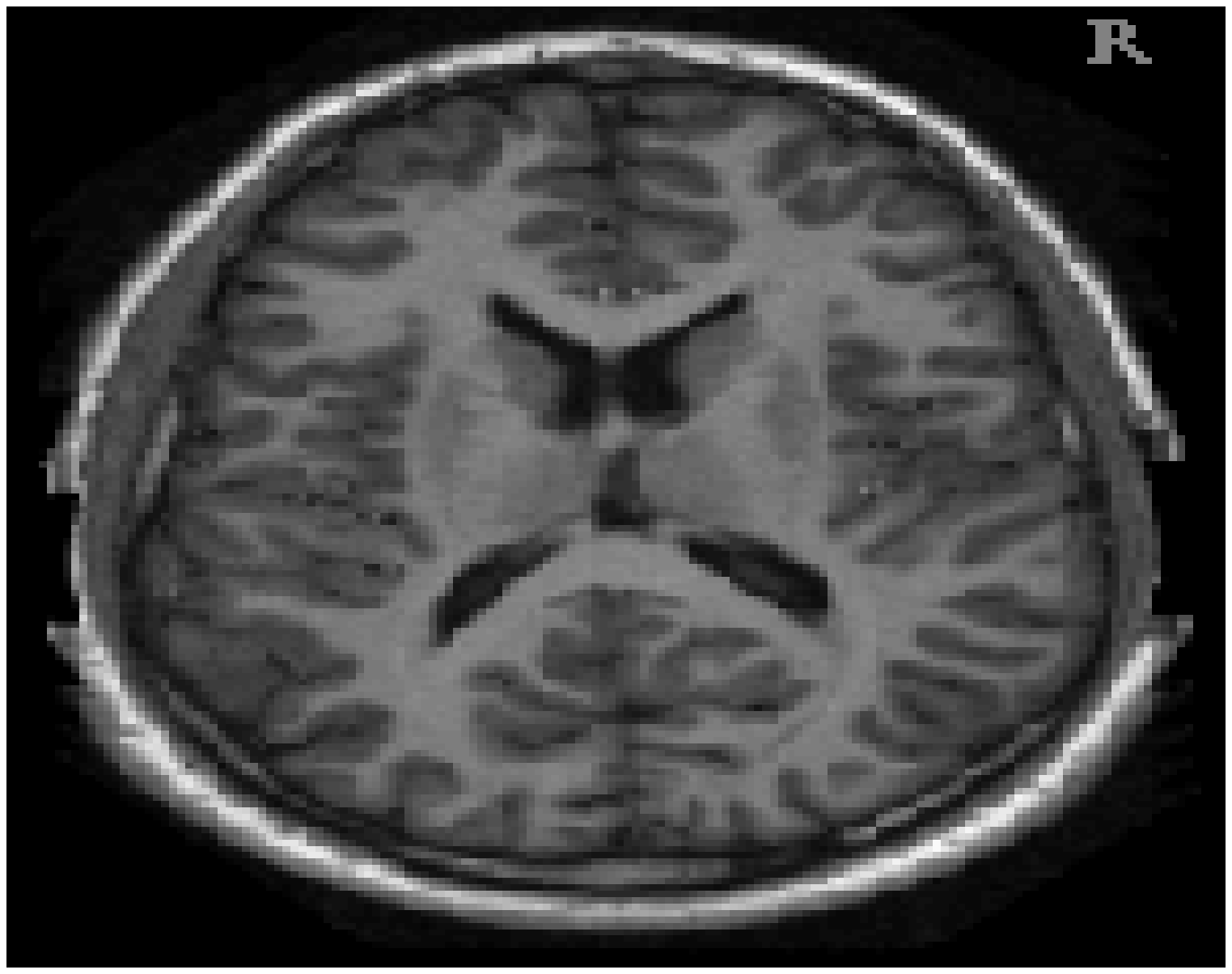}
\hspace{-0.3cm}
 \includegraphics[width=2cm,height=2cm]{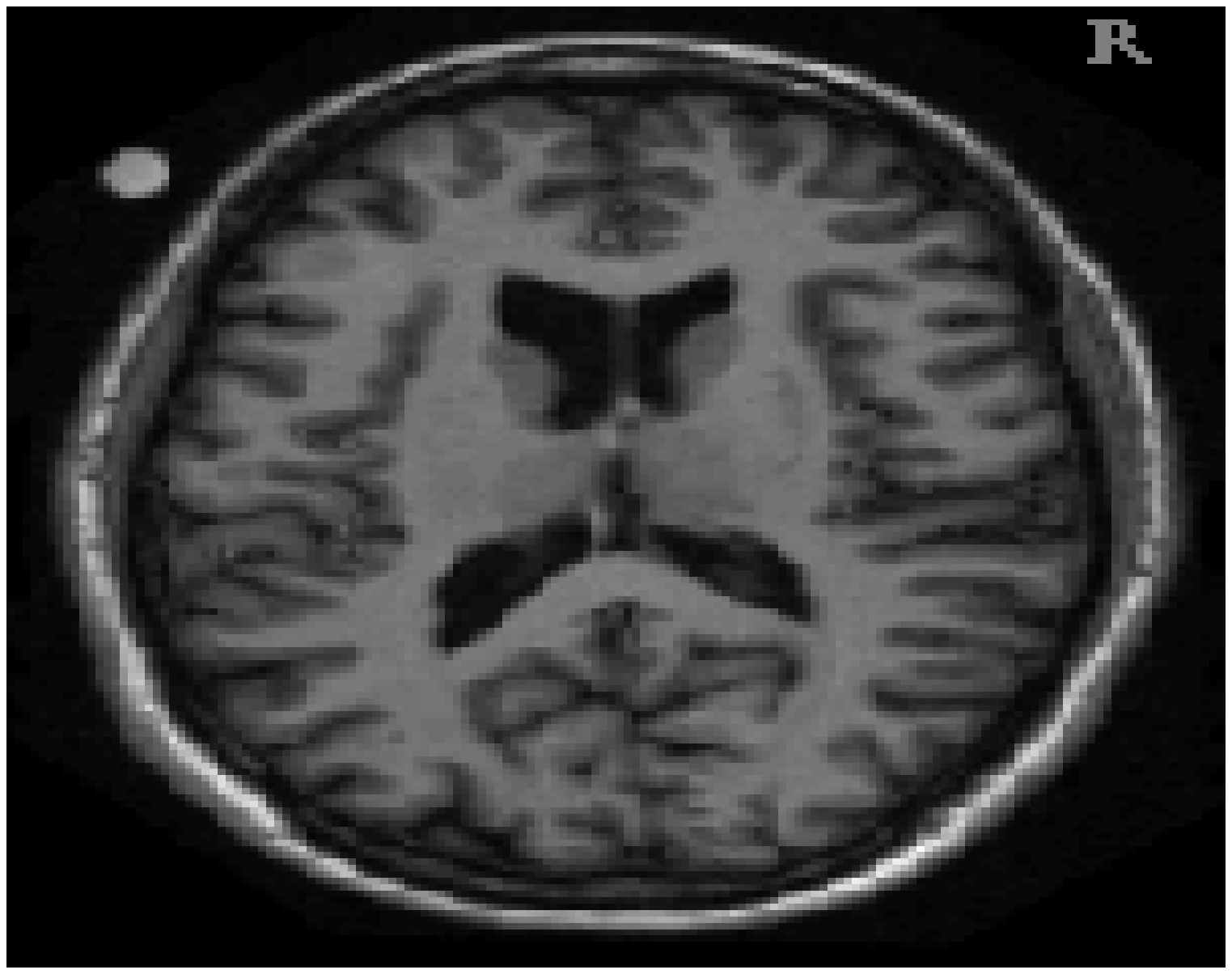}
\hspace{-0.3cm}
\includegraphics[width=2cm,height=2cm]{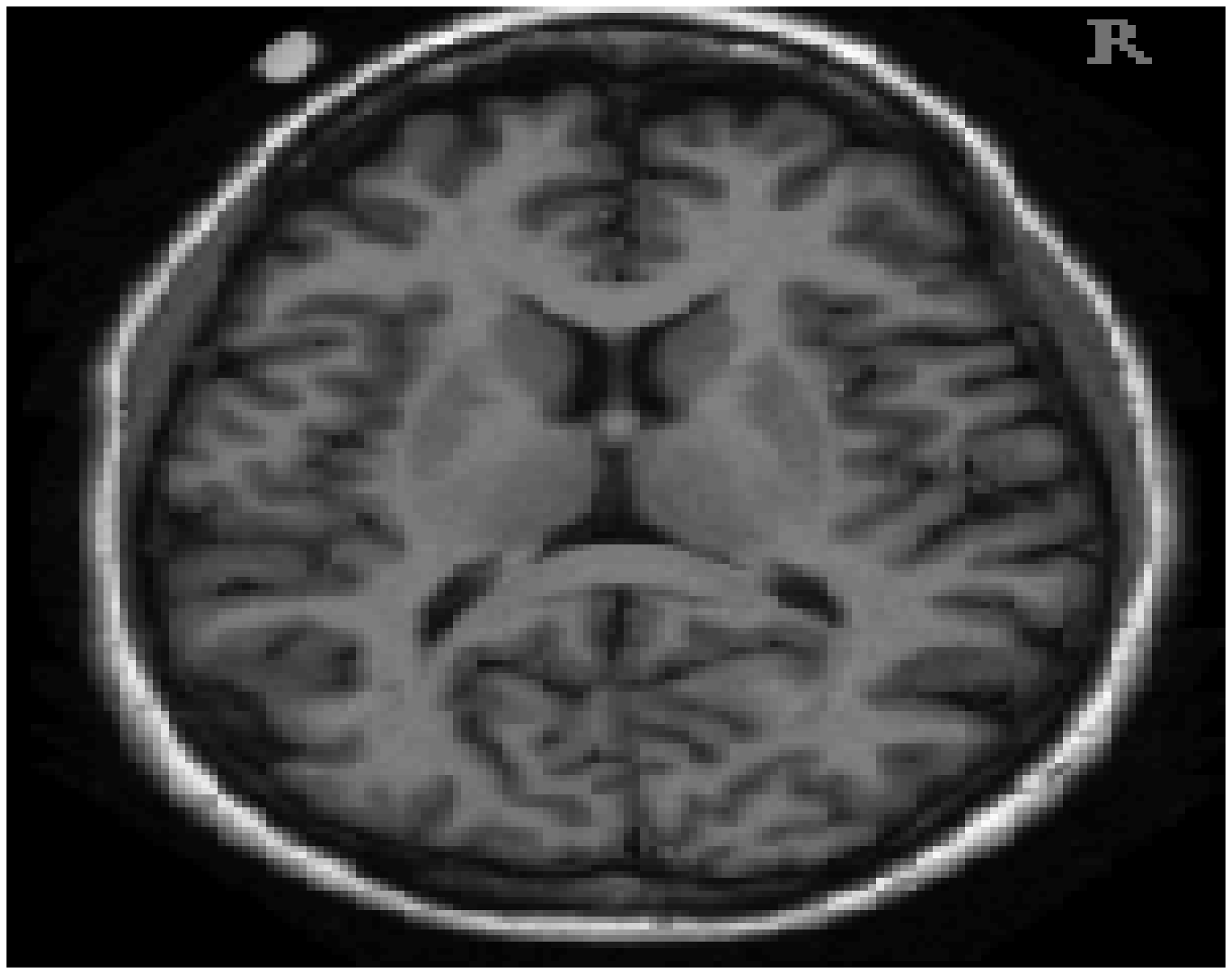}\\
\includegraphics[width=2cm,height=2cm]{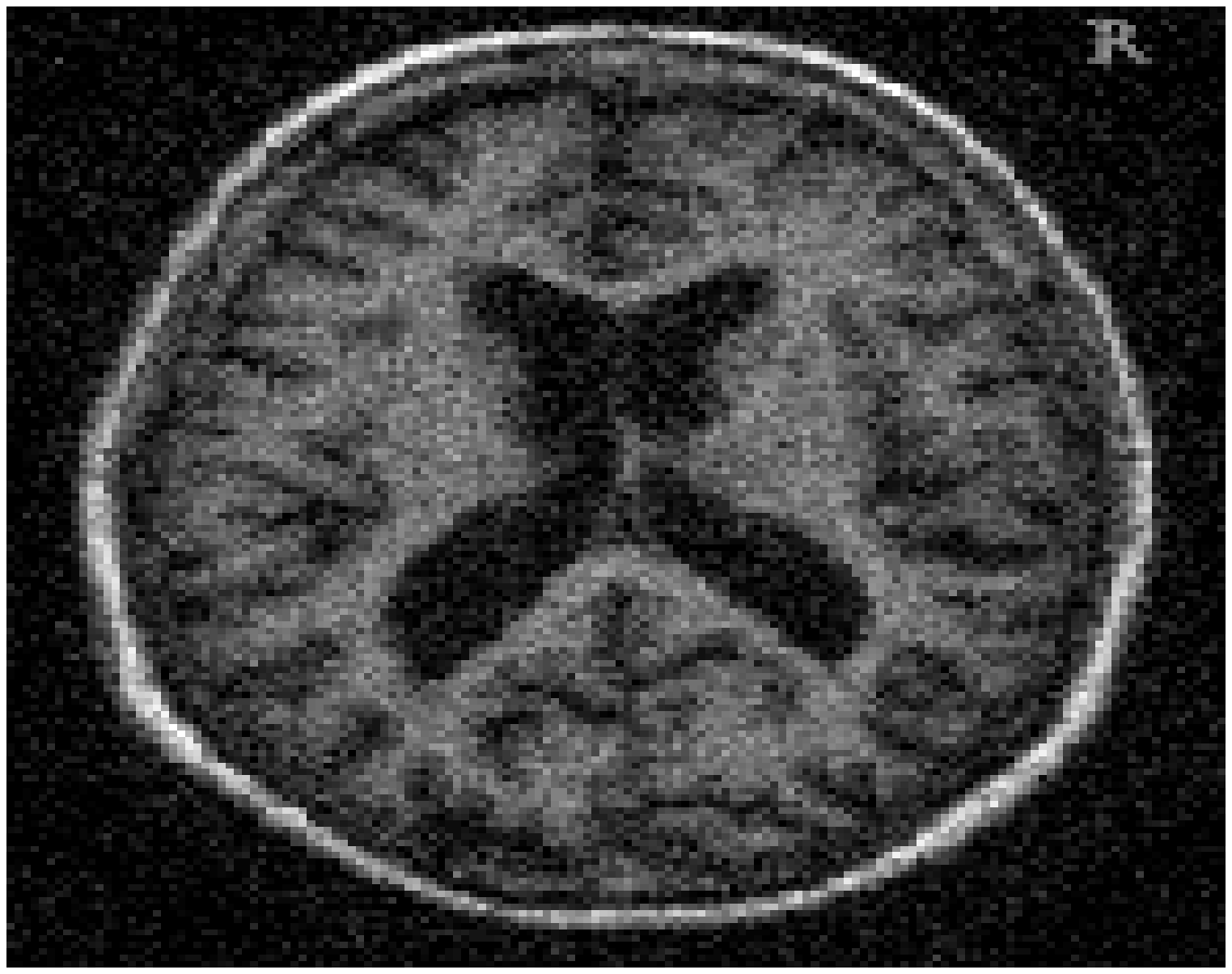}
\hspace{-0.3cm}
\includegraphics[width=2cm,height=2cm]{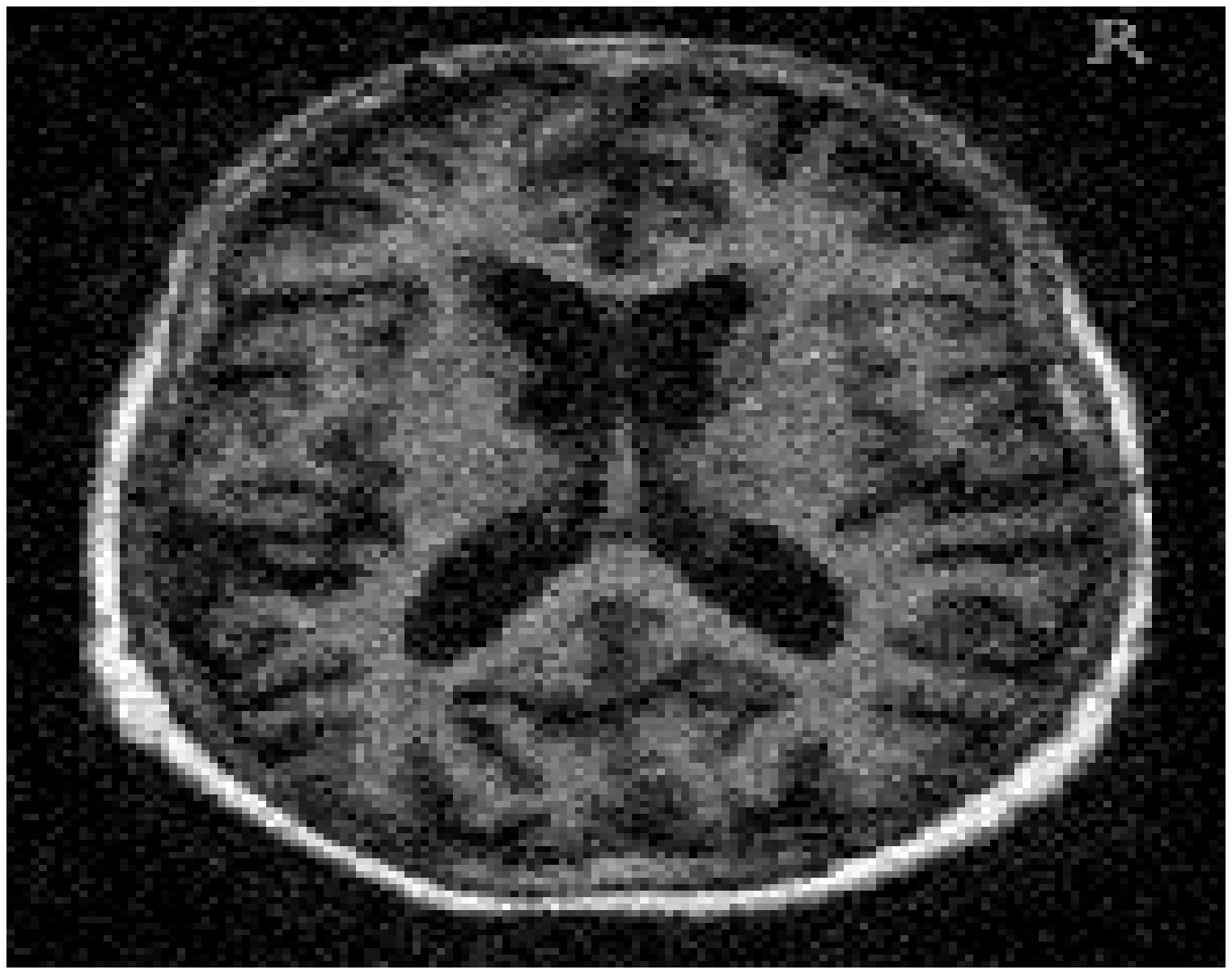}
\hspace{-0.3cm}
 \includegraphics[width=2cm,height=2cm]{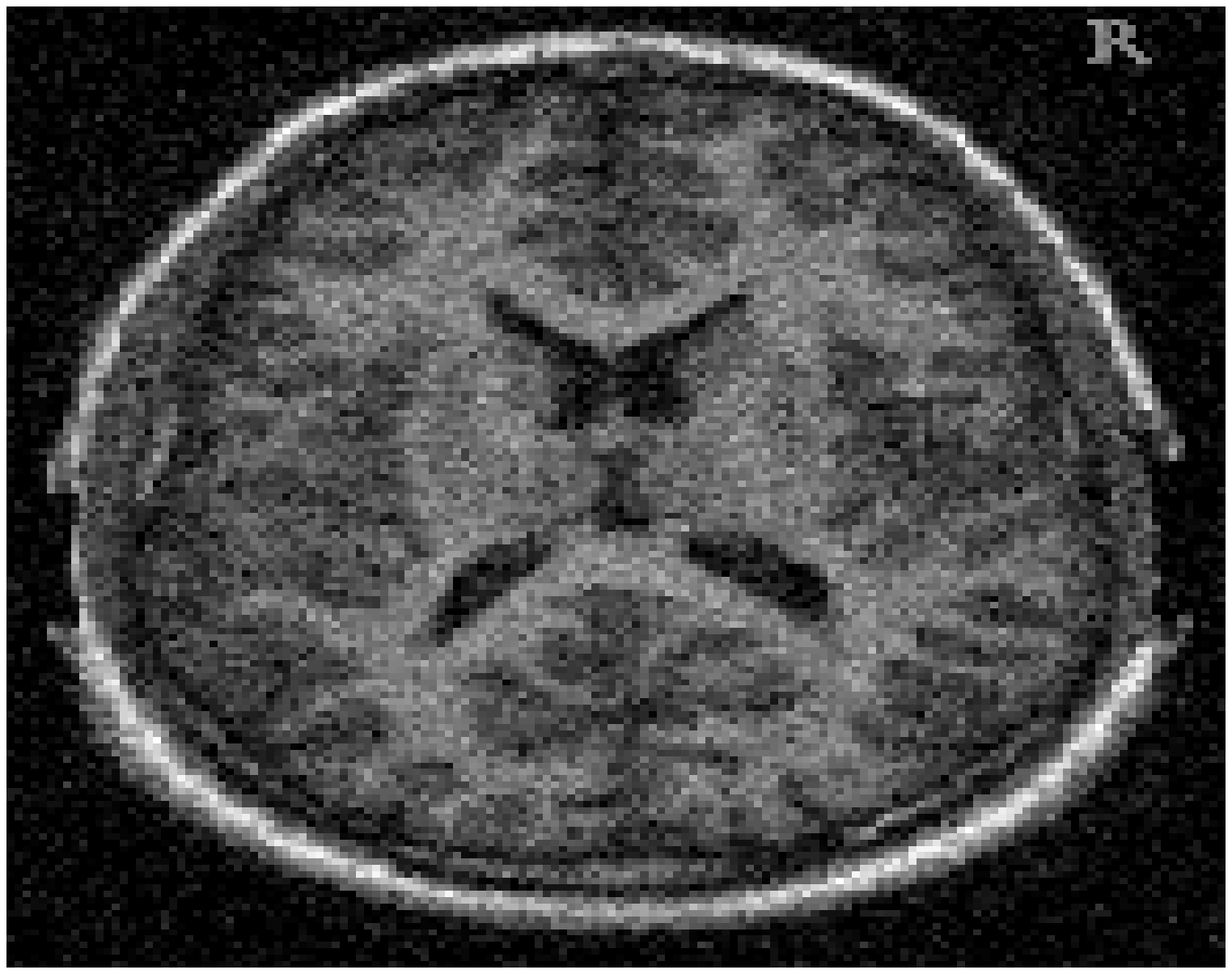}
\hspace{-0.3cm}
 \includegraphics[width=2cm,height=2cm]{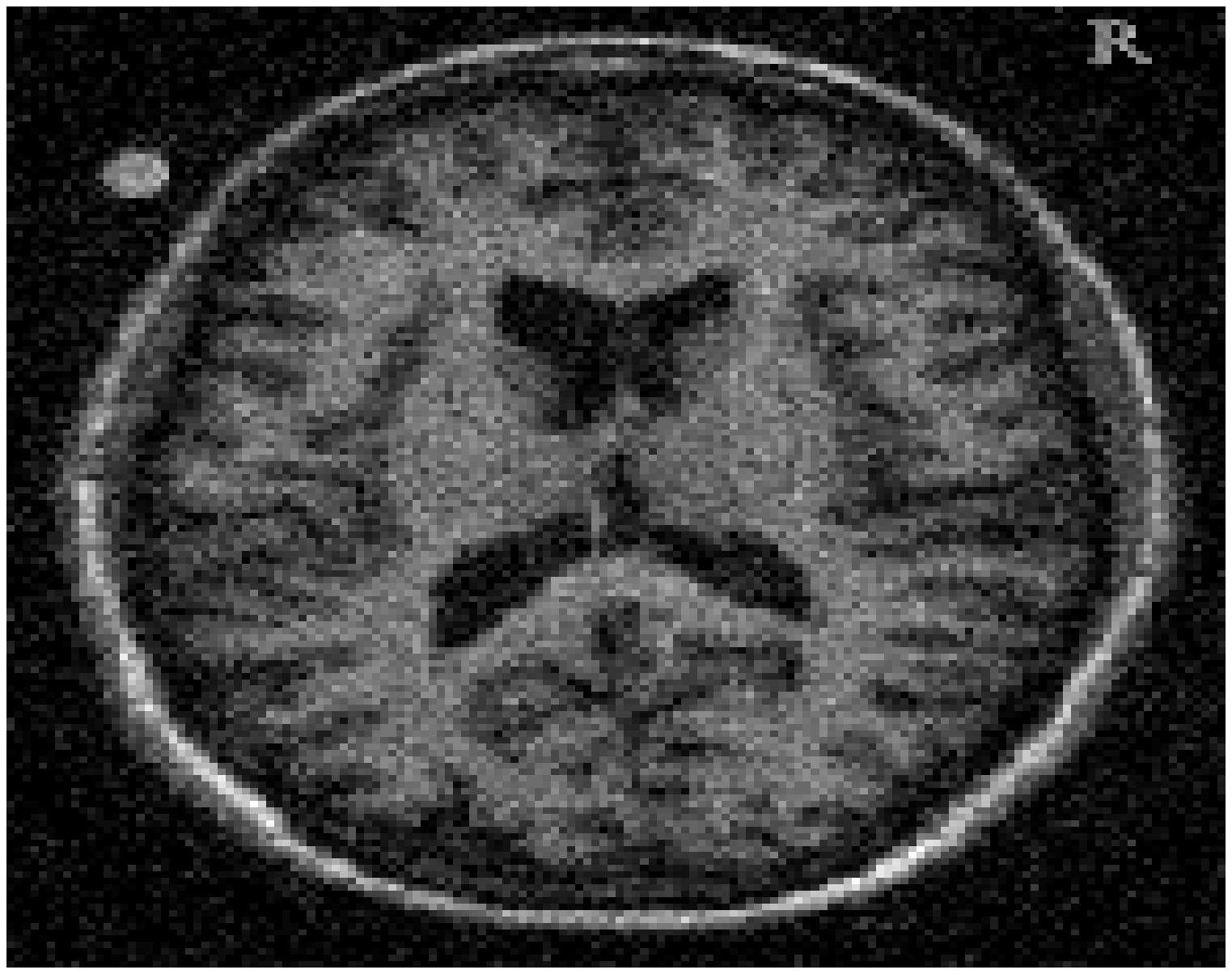}
\hspace{-0.3cm}
\includegraphics[width=2cm,height=2cm]{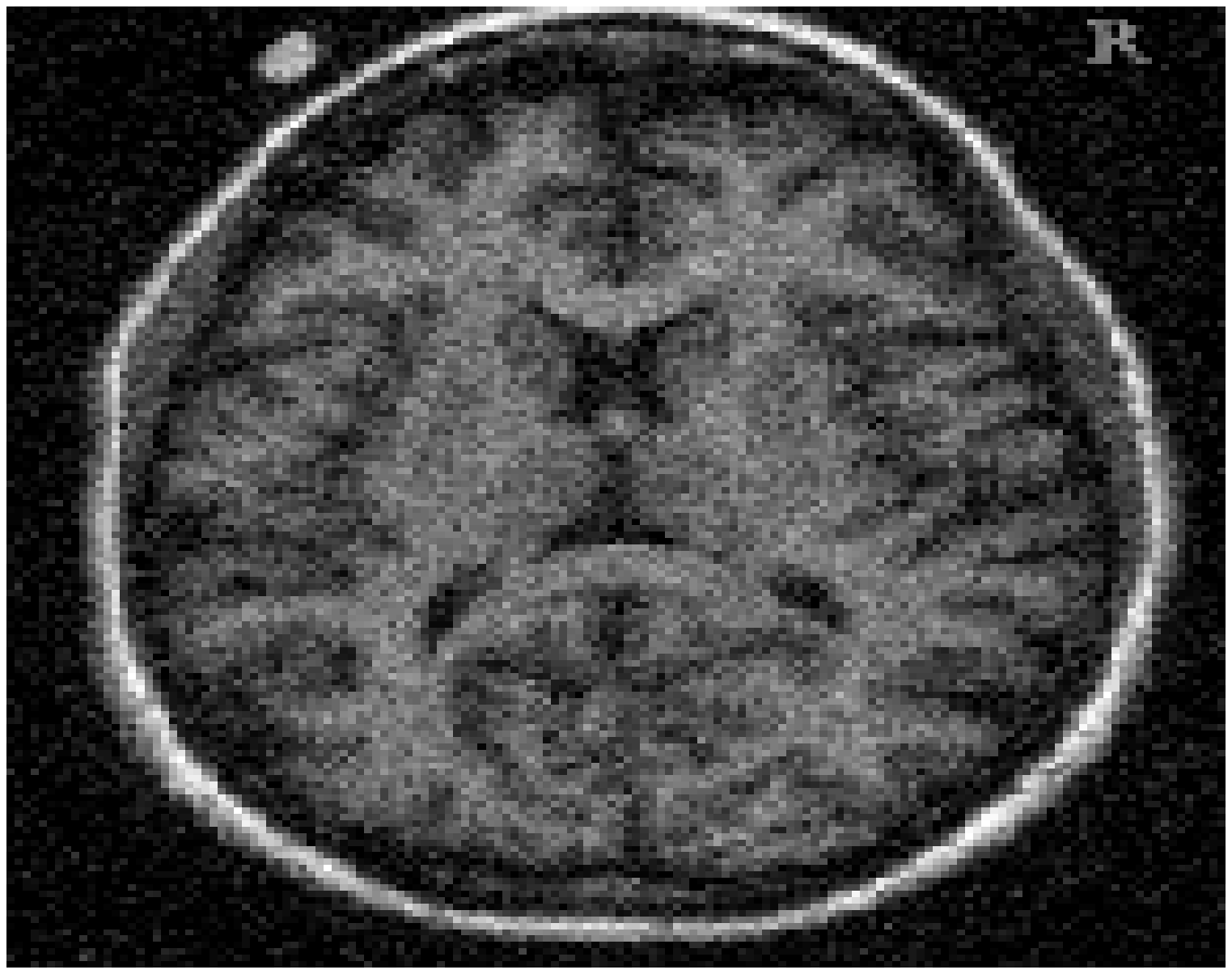}\\
\includegraphics[width=2cm,height=2cm]{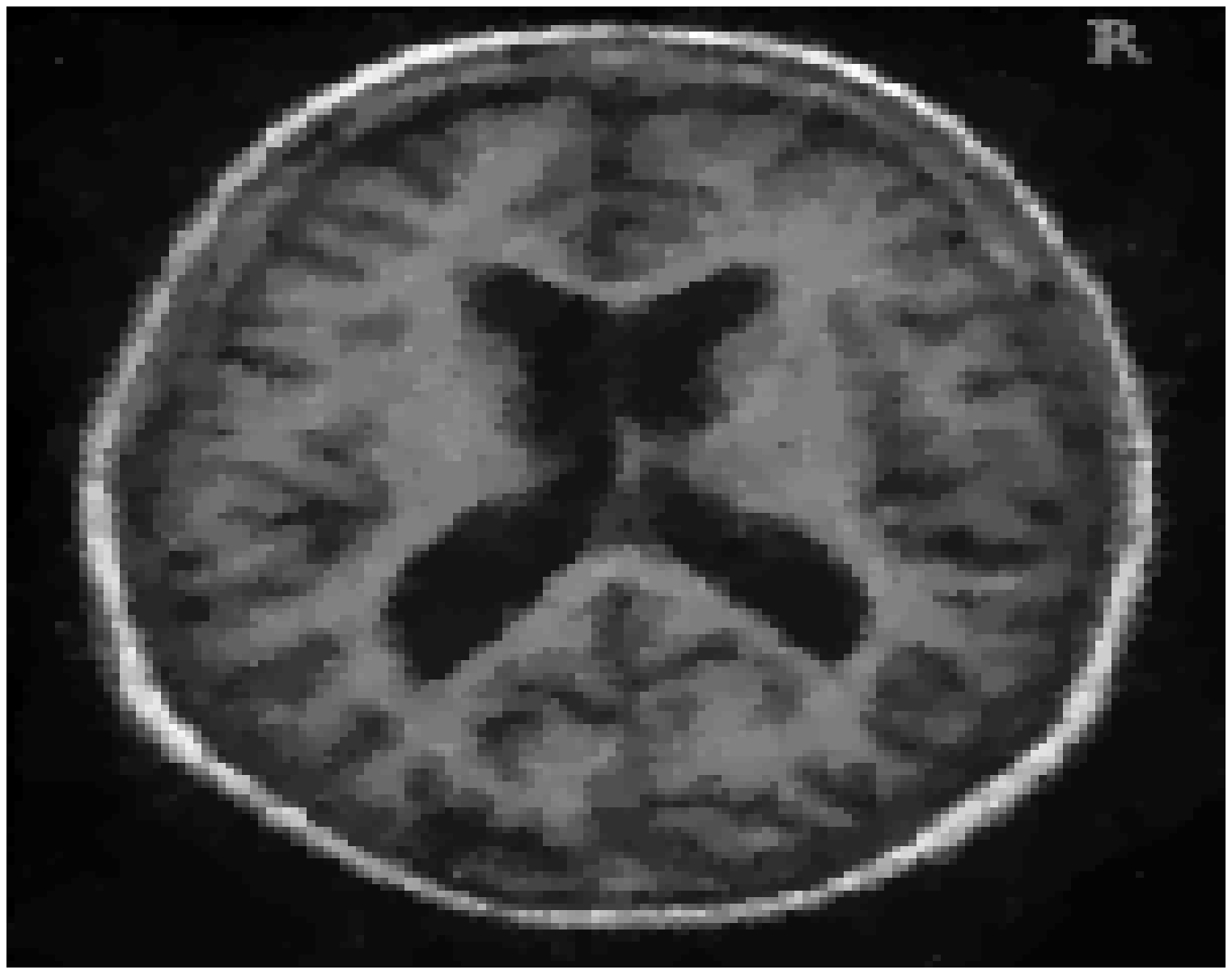}
\hspace{-0.3cm}
\includegraphics[width=2cm,height=2cm]{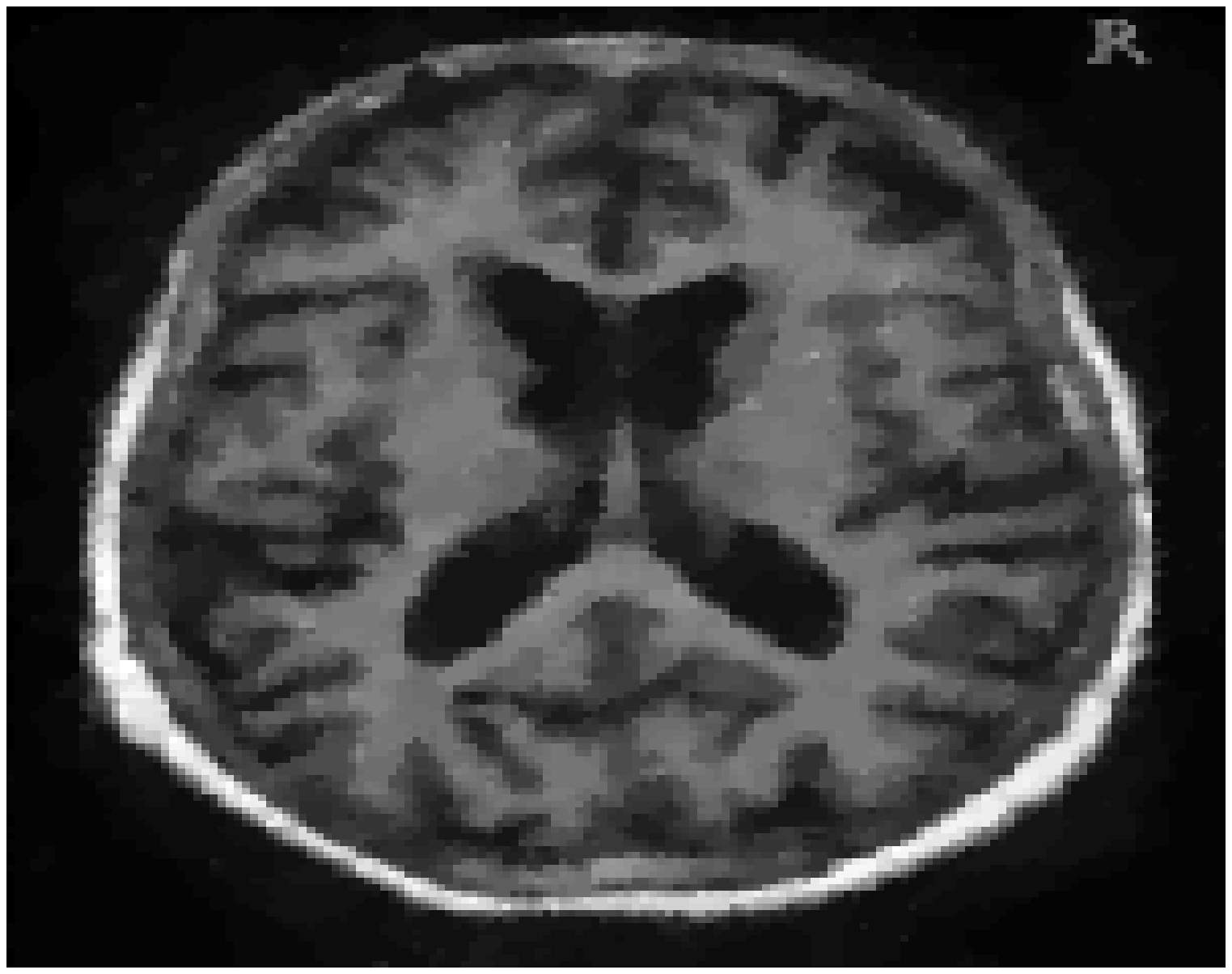}
\hspace{-0.3cm}
 \includegraphics[width=2cm,height=2cm]{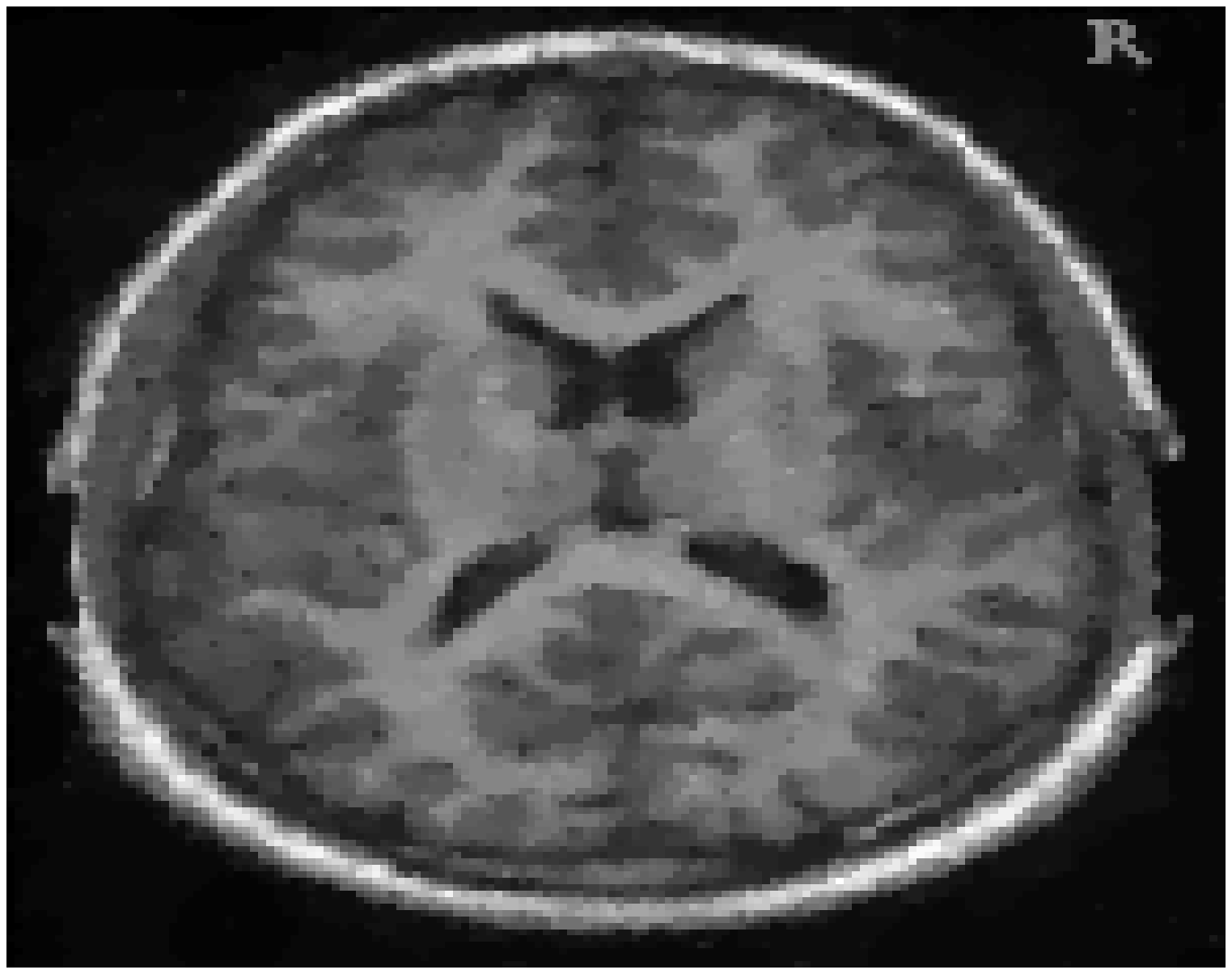}
\hspace{-0.3cm}
 \includegraphics[width=2cm,height=2cm]{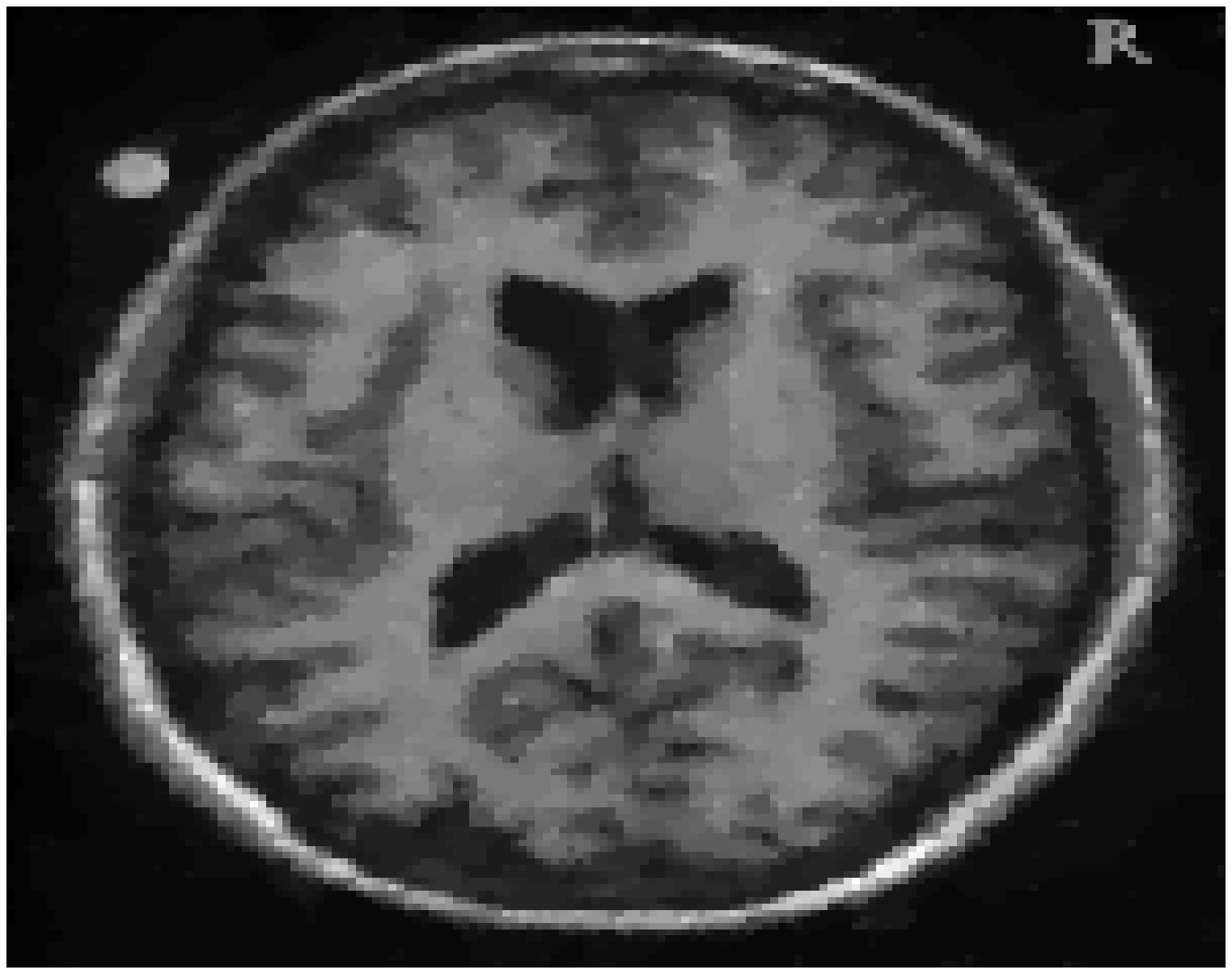}
\hspace{-0.3cm}
\includegraphics[width=2cm,height=2cm]{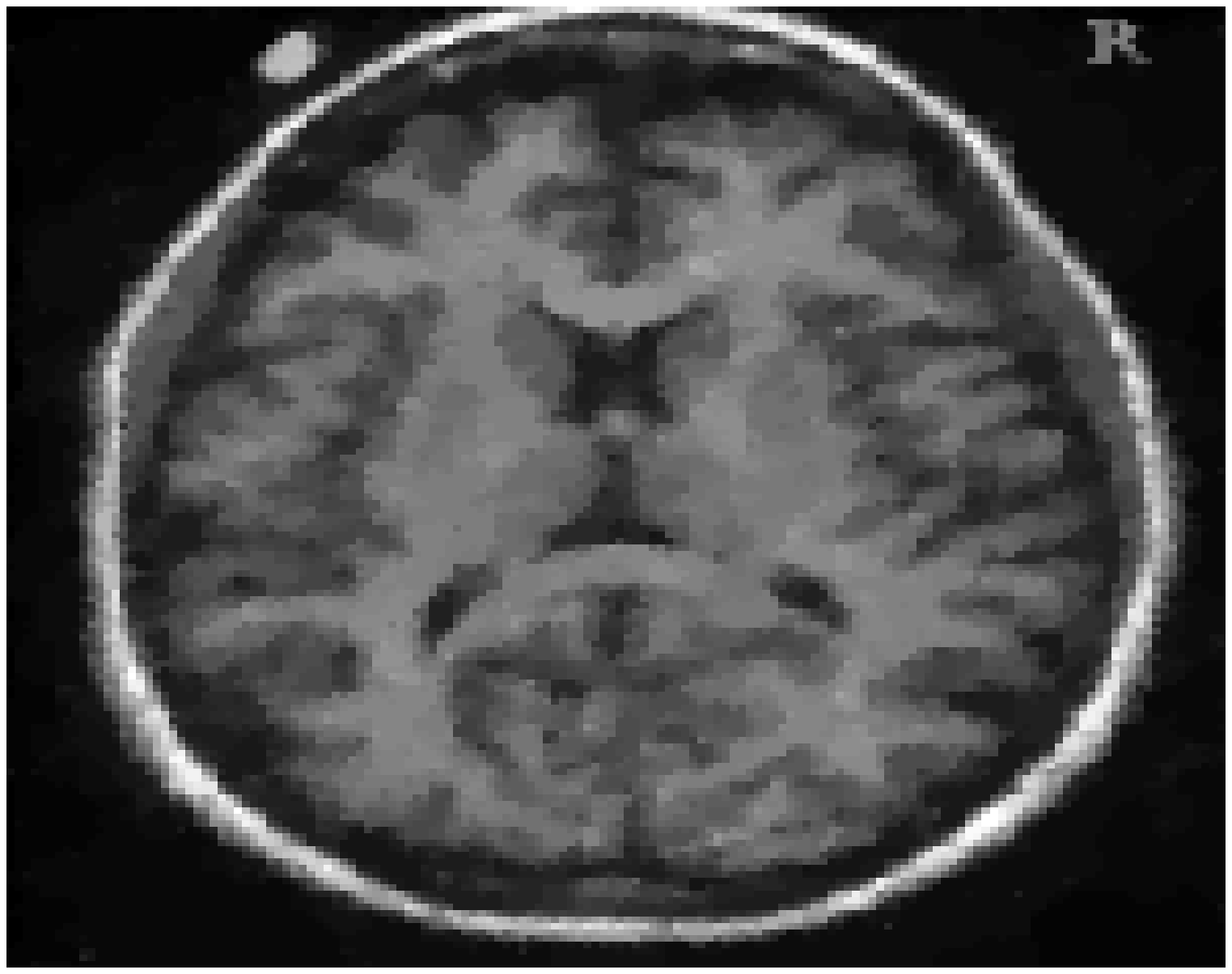}
\end{center}
\caption{Sample of $5$ images of a MRI brain database: original images (upper row), noisy images (middle row) and optimal denoised images (bottom row),  $\hat{\lambda_S}=3280.5$.}
\label{fig:training set}
\end{figure}

Figure \ref{fig:training set} shows an example of database of brain images\footnote{OASIS online database.} together with the optimal denoised version obtained by Algorithm 1 for single Gaussian noise estimation.

%\vspace{-1cm}
\paragraph{\textbf{Multiple noise estimation.}} We consider now a more interesting application of \eqref{optprob1}-\eqref{optprob2} where the image is corrupted by noises with different distributions. We consider the case where a combination of Gaussian and impulse noise is present. The fidelity term for the impulse distributed component is $\phi_1(u,f_k)=|u-f_k|$, whereas, as above, for the Gaussian noise we consider the fidelity $\phi_2(u,f_k)=(u-f_k)^2$, for every $k$. Each fidelity term is weighted by a parameter $\lambda_i, i=1,2$. Thus, we aim to solve:
\begin{equation}  \label{optprob1mixed}
\min_{(\lambda_1,\lambda_2),\ \lambda_i\geq 0}\,\frac{1}{2N}\sum_{k=1}^{N}\norm{\uu-u_k}^2_{L^2(\Omega)}
%+\frac{\beta}{2}\sum_{i=1}^{2} |\lambda_i|^2
\end{equation}
where, for each $k$, $\uu$ is now the solution of the regularised PDE:
\begin{equation} \label{optprob2mixed}
-\varepsilon\Delta \uu- \text{div}( h_\gamma(\nabla \uu)) +\lambda_1 h^1_\gamma (\uu-f_k)+\lambda_2 (\uu-f_k)=0,\quad k=1,\ldots,N. 
\end{equation}
In \eqref{optprob2mixed} the first and the second terms are as before while the third one corresponds to the Huber-type regularisation of $\text{sgn} (\uu-f_k)$. The adjoint state is computed as in \cite{delosreyes}, in a similar manner as \eqref{eq:optimality condition adjoint equation}. By taking also into account equations \eqref{batchgradient}-\eqref{batchgradientadj},
%\blue{Here I would do the same thing as for the Gaussian-noise constraint above.}
we solve \eqref{optprob1mixed}-\eqref{optprob2mixed} with $\varepsilon=10^{-12}, \gamma=100$  by means of Algorithm 1. 
%Note that in \eqref{optprob1mixed} the parameter $\beta$ is set equal to zero. 
%Again, in \eqref{optprob2mixed}, an elliptic-regularisation term is present. We follow faithfully \cite[Section 4.4]{delosreyes} and we get a modified version of the Euler-Lagrange equations of \eqref{optprob2mixed} that we solve by means of a semismooth Newton method.

We take as example slices of the brain database shown in Figure \ref{fig:training set} corrupted with both Gaussian noise distributed as $\mathcal{N}(0,0.005)$ and impulse noise with fraction of missing pixels $d=5\%$, and again solve \eqref{optprob1}-\eqref{optprob2} by solving the PDE constraints all at once and by using Dynamic Sampling for different $N$. In Table \ref{tablemixed} we report the results for the estimation of $\lambda_1$ and $\lambda_2$.

%\begin{figure}[!h]
%\begin{subfigure}{0.45\textwidth}
%\centering
%\includegraphics[height=3cm,width=5cm]{samplesizemixed.eps}
%\caption{\quad Dynamic sample size increasing.}
%\end{subfigure}
%\hspace{0.5cm}
%\begin{subfigure}{0.45\textwidth}
%\centering
%\vspace{0.5cm}
%\includegraphics[height=2.5cm,width=1.7cm]{sliceor.png}
%\includegraphics[height=2.5cm,width=1.7cm]{slicen.png}
%\includegraphics[height=2.5cm,width=1.7cm]{sliced.png}
%\vspace{0.2cm}
%\caption{Example of database and optimal denoised brain image.}
%\end{subfigure}
%\caption{Optimal parameters selection and optimal denoising result for \eqref{optprob1}-\eqref{optprob2}. $|S_0|=20\% N,\ \theta=0.3$. }
%\label{sizesamplemixed}
%\end{figure}

\begin{center}
{\scriptsize
    \begin{tabular}{| c | c | c | c | c | c | c | c | } 
    \hline
$N$ & $\hat{\lambda_1}_S$ & $\hat{\lambda_2}_S$ & $|S_0|$ & $|S_{end}|$ & eff. & eff. Dyn.S. &  diff. \\ \hline
$10$ & $86.31$ & $28.43$ & $2$ & $7$ & $180$ & \bm{$70$} & $5.2\%$  \\ \hline
$20$ &  $90.61$ & $26.96$ & $4$ & $6$ & $920$ & \bm{$180$}  & $5.3\%$ \\ \hline
$30$ & $94.36$ & $29.04$  & $6$ & $7$ & $2100$  & \bm{$314$} & $5.6\%$ \\ \hline
$40$ & $88.88$ & $31.56$  & $8$  & $8$ & $880$ &  \bm{$496$}  & $1.2\%$ \\ \hline
$50$ & $88.92$  & $29.81$ & $10$ & $10$ & $2200$ & \bm{$560$}  & $<1\%$  \\ \hline
$60$ & $89.64$  & $28.36$ & $12$ & $12$ & $1920$ & \bm{$336$} & $1.9\%$  \\ \hline
$70$ & $86.09$  & $28.09$ & $14$ & $14$ & $2940$ & \bm{$532$}  & $3.3\%$  \\ \hline
$80$ & $87.68$  & $29.97$ & $16$ & $16$ & $3520$ & \bm{$448$}  & $<1\%$  \\ \hline
 \end{tabular}
 }
\captionof{table}{$\hat{\lambda_1}_S$ and $\hat{\lambda_2}_S$ are the optimal weights for \eqref{optprob1mixed}-\eqref{optprob2mixed} estimated with Dynamic Sampling. We observe again a clear improvement in efficiency (i.e. number of PDEs solved). As above, $|S_0|=20\% N$ and $\theta=0.5$. }
\label{tablemixed}
\end{center}

\begin{figure}[!h]
\begin{center}
\includegraphics[height=3.2cm]{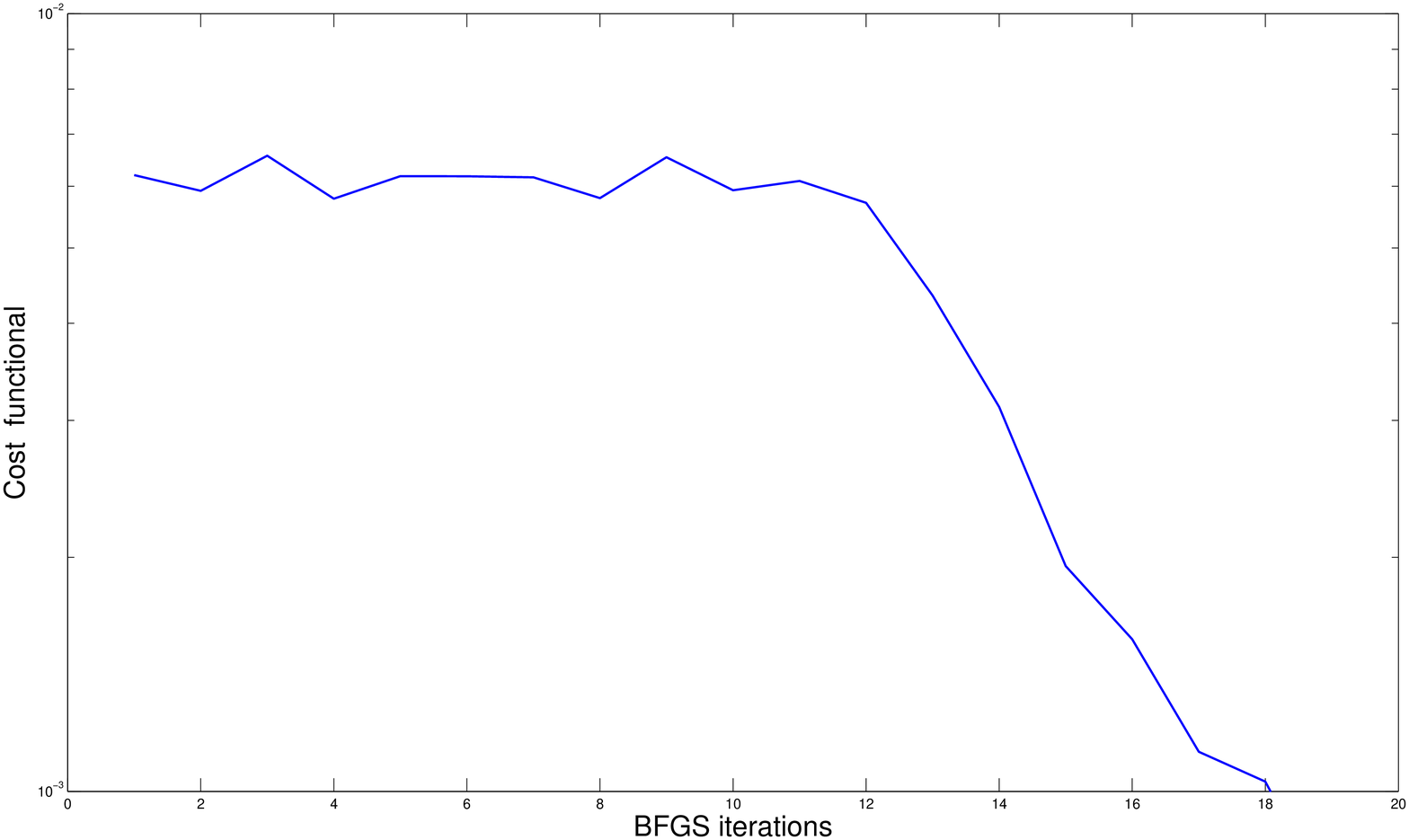}
\includegraphics[height=3.2cm]{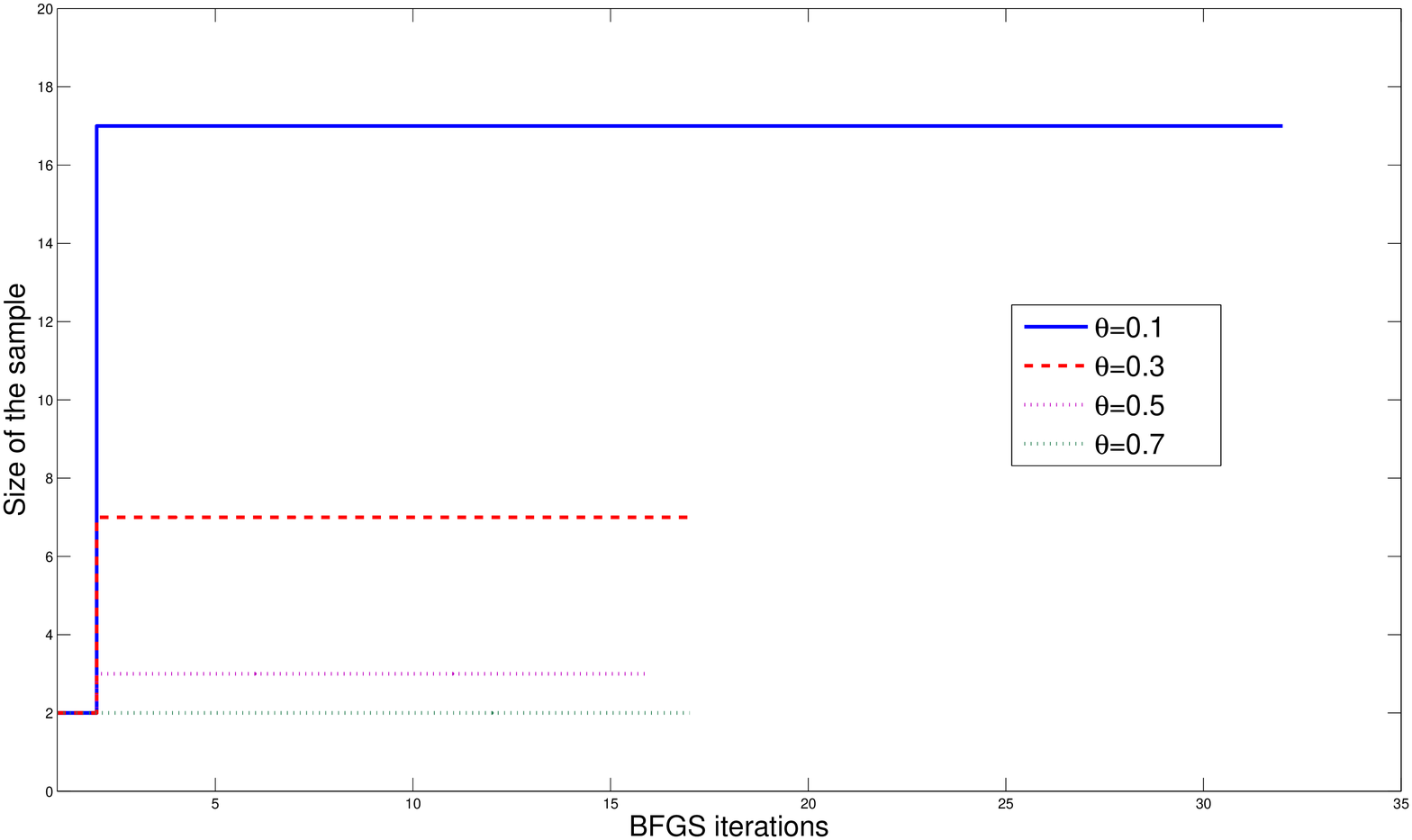}
\end{center}
\caption{\emph{Left}: evolution of BFGS with Dynamic Sampling along the iterations. \emph{Right}: samples size changes in Algorithm 1 for different values of $\theta$. For each value of $\theta$, the result is plotted till convergence. For this example $N=20, |S_0|=2$.}
\label{graphs}
\end{figure}

\paragraph{\textbf{Convergence and sensitivity.}} Figure \ref{graphs} shows two  features of Algorithm 1 applied to solve problem \eqref{optprob1gauss}-\eqref{optprob2gaussreg}. On the left we represent the evolution of the cost functional along the BFGS iterations. Because of the sampling strategy, in the early iterations of BFGS the problem considered varies quite a lot, thus showing oscillations. Once evolving the process, the convergence is superlinear. On the right we represent the sensitivity with respect to the accuracy parameter $\theta$ (cf. \eqref{descentdirectionvar}): smaller values of $\theta$ penalise larger variances on $\nabla J_S$, thus favouring larger samples. Larger values of $\theta$  allow larger variances on $\nabla J_S$ and, consequently, smaller sample sizes. In this case, efficiency improves, but accuracy might suffer as shown in Table \ref{tableaccur}.
%\vspace{-0.5cm}

\begin{center}
{\scriptsize
    \begin{tabular}{| c | c| c |} 
    \hline
 $\quad \theta\quad $ \quad & efficiency \quad &\quad difference \quad \\ \hline
$0.1$ & $516$ & $0.07\%$  \\ \hline
$0.3$ &  $246$ & $4.3\%$  \\ \hline
$0.5$ &  $92$ &  $5.9\%$\\ \hline
$0.7$ &   $68$ &  $15\%$\\ \hline

    \end{tabular}}
\captionof{table}{As $\theta$ increases we observe  improvements upon the efficiency as smaller samples are allowed. However, the relative difference with the value estimated without sampling shows that accuracy suffers.}
\label{tableaccur}
\end{center}

\paragraph{Conclusions.} In this paper, we propose an efficient and competitive technique to solve numerically the constrained optimisation problem \eqref{optprob1}-\eqref{optprob2} designed for learning the noise model in a TV denoising framework accounting for different types of noise. The set of nonsmooth PDE constraints resembles a large-size training database of clean and noisy images that allows a more robust estimation of parameters. To solve the problem, we use \emph{Dynamic Sampling} methods, proposed in \cite{byrd} for \emph{linear} constrained problems. The idea consists in selecting just a small sample of the PDEs that need to be solved over the whole database and then, during the progression of the algorithm, verify whether such a size produces approximations that are accurate enough. Extended to our nonlinear framework, the results show a remarkable improvement in  efficiency, which reflects in reduced computational times for both single noise estimations as well as for mixed ones. Further directions for future research are an accurate analysis of convergence properties of such a scheme as well as the design of a similar algorithm  for the case of a $L^1$-regularisation on the parameter vector.

\paragraph{Acknowledgements.}
This work is supported by the King Abdullah University for Science and Technology (KAUST) Award No. KUK-I1-007-43, the EPSRC first grant Nr. EP/J009539/1 and the Cambridge Center for Analysis (CCA).
%Carola-Bibiane Sch\"{o}nlieb acknowledges financial support provided by the Cambridge Centre for Analysis (CCA), the Royal Society International Exchanges Award IE110314 for the project \emph{High-order Compressed Sensing for Medical Imaging}, the EPSRC first grant Nr. EP/J009539/1 \emph{Sparse \& Higher-order Image Restoration}, and the EPSRC / Isaac Newton Trust Small Grant on \emph{Non-smooth geometric reconstruction for high resolution MRI imaging of fluid transport in bed reactors}. Further, this publication is based on work supported by Award No.
%KUK-I1-007-43, made by King Abdullah University of Science and Technology (KAUST).

\end{document}